\documentstyle{article}

\catcode`\@=11
\font\tenmsx=msxm10
\font\sevenmsx=msxm7
\font\fivemsx=msxm5
\font\tenmsy=msym10
\font\sevenmsy=msym7
\font\fivemsy=msym5
\newfam\msxfam
\newfam\msyfam
\textfont\msxfam=\tenmsx  \scriptfont\msxfam=\sevenmsx
\scriptscriptfont\msxfam=\fivemsx   
\textfont\msyfam=\tenmsy  \scriptfont\msyfam=\sevenmsy
\scriptscriptfont\msyfam=\fivemsy
\def\hexnumber@#1{\ifnum#1<10 \number#1\else
\ifnum#1=10 A\else\ifnum#1=11 B\else\ifnum#1=12 C\else
\ifnum#1=13 D\else\ifnum#1=14 E\else\ifnum#1=15 F\fi\fi\fi\fi\fi\fi\fi}
\def\msx@{\hexnumber@\msxfam}
\def\msy@{\hexnumber@\msyfam}
\mathchardef\nsubseteq="3\msy@2A
\mathchardef\nsubseteqq="3\msy@22
\catcode`\@=\active

\newtheorem{theorem}{Theorem}[section]
\newtheorem{defn}[theorem]{Definition}
\newtheorem{prop}[theorem]{Proposition}
\newtheorem{lem}[theorem]{Lemma}

\newtheorem{eple}[theorem]{Example}

\newcommand{\implies}{\rightarrow}
\newcommand{\sub}{\subseteq}
\newcommand{\sups}{\supseteq}

\newcommand{\sm}{\setminus}
\newcommand{\proof}{{\bf Proof:}\ }
\newcommand{\qed}{\hfill $\Box$}
\newcommand{\noi}{\noindent}

\newcommand{\ra}{\rangle}
\newcommand{\la}{\langle}

\newcommand{\rb}{\mbox{{\large )}}}
\newcommand{\lsb}{\mbox{{\large [}}}

\newcommand{\CA}{{\cal A}}
\newcommand{\CB}{{\cal B}}
\newcommand{\CC}{{\cal C}}

\newcommand{\CF}{{\cal F}}
\newcommand{\CG}{{\cal G}}
\newcommand{\CH}{{\cal H}}
\newcommand{\CI}{{\cal I}}

\newcommand{\CM}{{\cal M}}
\newcommand{\CO}{{\cal O}}

\newcommand{\CS}{{\cal S}}  
\newcommand{\CT}{{\cal T}}   
\newcommand{\CU}{{\cal U}}  
 
\newcommand{\CX}{{\cal X}} 
  
\newfont{\gothic}{eufm10} 
\newcommand{\gb}{{\gothic b }}
\newcommand{\gc}{{\gothic c }} 
\newcommand{\gd}{{\gothic d }}
\newcommand{\fo}{\; ^\omega  \! \omega}
\newcommand{\fso}{[\omega]^{<\omega}}
\newcommand{\ff}{\; ^{\omega} \! 2}
\newcommand{\fff}{\; ^{<\omega} \! 2}
\newcommand{\so}{[\omega]^\omega}
\newcommand{\po}{\mbox{{\Large $\wp$}}(\omega)}
\newcommand{\fa}{\forall}
\newcommand{\ex}{\exists}
\newcommand{\fai}{\forall^\infty}

\title{Strong meager properties for filters}
	 \date{}
	 \author{Claude Laflamme\thanks{This research was partially supported by
         NSERC of Canada.  \newline
         1980 {\em AMS Subj. Class. (1985 Revision)} Primary 04A20; Secondary 
         03E05,03E15,03E35.
         } \\
         Department of Mathematics and Statistics \\
          University of Calgary \\
         Calgary, Alberta  \\
         Canada T2N 1N4}
						  
\begin{document}
\maketitle
\begin{abstract}
     We analyze several ``strong meager'' properties for filters on the
natural numbers between the classical Baire property and a filter being
$F_\sigma$.  Two such properties have been studied by Talagrand and a
few more combinatorial ones are investigated.  In particular, we define
the notion of a P$^+$-filter, a generalization of the traditional
concept of P-filter, and prove the existence of a non-meager
P$^+$-filter.  Our motivation lies in understanding the structure of
filters generated by complements of members of a maximal almost disjoint
family. 
\end{abstract}
							  
\section{Introduction}

 We investigate several ``smallness'' properties of filters as
strengthening of the meager property.  Indeed, although a meager filter
might be thought of as small, it may happen for example that adding one
set to a meager filter yields an ultrafilter (see the section on
examples).  Our motivation lies in the unresolved relationship between
the minimum size of a maximal almost disjoint (mad) family and the
minimum size of a dominating family; the relation with the minimum size
of an unbounded family has been settled by Shelah in \cite{shelah}. 
Part of the problem lies in the technical question whether filters
generated by the complement of members of a mad family can be
diagonalized in a forcing extension without adding an unbounded real. 
This loosely means that we must understand the different possibilities
for a filter to code an unbounded function when diagonalized.  Very
partial results were obtained in \cite{laf} where it became apparent
that strong meager properties are required on filters for this program
to succeed.  It is the goal of this paper to analyze some of these
stronger meager properties.

Our terminology is standard but we review the main concepts and
notation.  The natural numbers will be denoted by $\omega$, $\ff$ and
$\fo$ denote the collection of functions from $\omega$ to 2 and to $\omega$
respectively; similarly, $\po$ and $\so$ denote the collection of all
and infinite subsets respectively.  Two orderings are standard, the
first one for subsets of $\omega$ is $A \sub^*B$ if $A \sm B$ is finite,
the second one for functions is $f\leq^*g$ if $f(n) \leq g(n)$ for all
but finitely many $n$.  We can view members of $\po$ as members of $\ff$
by considering their characteristic functions.  We equip $\ff$ with the
product toplogy of the discrete space $\{0,1\}$.  A basic neighbourhood
is then given by sets of the form

\[ \CO_s = \{f \in \ff: s \sub f \} \] 

\noi where $s \in \fff$, the collection of finite binary sequences.  The
terms ``nowhere dense'', ``meager'', ``Baire property'' all refer to
this topology.

A filter is a collection of subsets of $\omega$ closed under finite
intersections, supersets and containing all cofinite sets; it is called
proper if it contains only infinite sets.  Given a collection $\CX \sub
\po$, we let $\la \CX \ra$ denote the filter generated by $\CX$.  For a
filter $\CF$, $\CF^+$ denotes the collection of all sets $X$ such that
$\la \CF,X \ra$ is a proper filter; it is useful to notice that $X \in
\CF^+$ if and only if $X^c \notin \CF$.  We say that $X \in \so$
diagonalizes a filter $\CF$ if $X \sub^* Y$ for all $Y \in \CF$.  The
{\em Fr\'echet} filter is the collection of cofinite sets, denoted by
{\sl Fr}.  Given two filters $\CF$ and $\CG$, we form a third filter by

\[ \CF \otimes \CG = \{X \sub \omega \times \omega:
           \{n: \{m: (n,m) \in X\} \in \CG \} \in \CF \} \}. \]

\noi This filter can be viewed as a filter on $\omega$ if desired by fixing a
bijection between $\omega \times \omega$ and $\omega$. 

The following important result characterizes meager filters in terms
of combinatorial properties.

\begin{prop}
(Talagrand (\cite{tal1})) The following are equivalent for a filter $\CF$:
\begin{enumerate}
\item $\CF$ has the Baire property.
\item $\CF$ is meager.
\item There is a sequence $n_0 < n_1 < \cdots$ such that
   \[ (\fa X \in \CF)(\fai k) \; X \cap [n_k,n_{k+1}) \neq \emptyset. \]
\end{enumerate}
\end{prop}

A family $\CF \sub \fo$ is called bounded if
\[ (\ex g \in \fo)(\fa f \in \CF) \; f \leq^* g. \]
\noi Similarly, a family $\CF \sub \fo$ is called dominating if
\[ (\fa g \in \fo)(\ex f \in \CF) \; g \leq^* f. \]

\noi The cardinals \gb and \gd denote respectively the minimum cardinality
of an unbounded or dominating family.

 \section{Strong meager properties}

\begin{defn}
Let $\CF$  be a filter on $\omega$.
\begin{enumerate}

\item $\CF$ is called {\it completely meager} (Talagrand, \cite{tal2})if the filter
$\la \CF,X \ra$ is meager whenever $X^c \notin \CF$. 

\item $\CF$ is called {\it hereditarily meager} if for each $f \in \fo$, the 
filter $f(\CF)$ is meager (or improper).

\item $\CF$ is called {\it weakly hereditarily meager}  if for each
$f \in \fo$, either $f^{-1}(i)$ is compatible with $\CF$ for some $i$, or else
$f(\CF)$ is meager.

\item $\CF$ has the {\it strong Baire property} (Talagrand, \cite{tal2}) if $\CF
\cap \CC$ has the Baire property relative to $\CC$ for each closed $\CC \sub \po$.

\item $\CF$ is a P$^+$-filter if given any decreasing sequence $\la X_0 \;
^*\!\sups X_1\; ^*\!\sups X_2 \; ^*\!\sups \cdots \ra$ from $\CF^+$, there is
an $X \in \CF^+$ such that $X \sub^* X_n$ for each $n$. 

\item $\CF$ is a strong P$^+$-filter if given any sequence $\la X_0,
X_1, X_2 ,\cdots \ra$ from $\CF^+$, there is an $f \in \fo$ such that for
all $X \in \CF$, $X\cap X_n \cap f(n) \neq \emptyset$ for all but
finitely many $n$. 

\end{enumerate}
\end{defn}

We shall prove the following implications between these notions.  In a
subsequent section, we shall show by examples that many of these 
implications do not reverse. Some open questions remain.

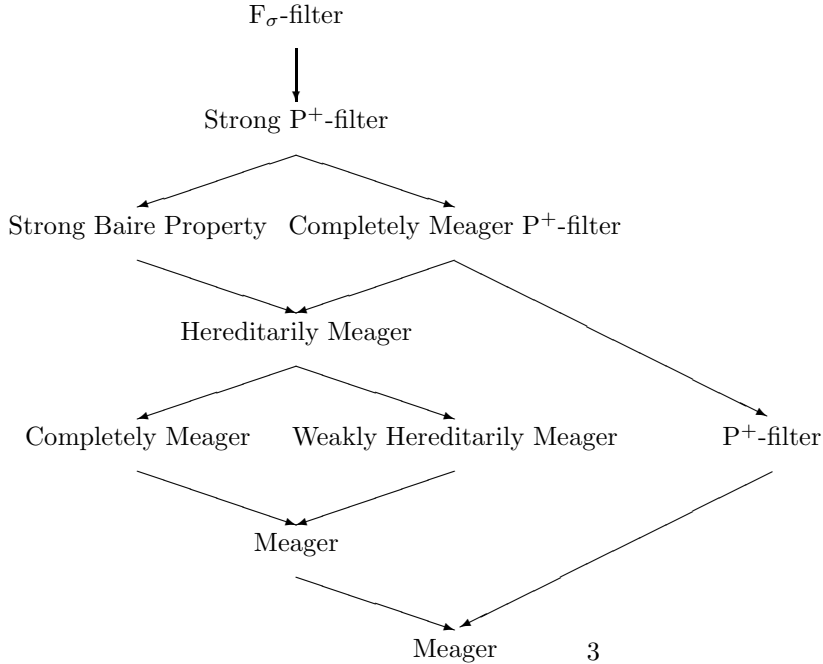
\begin{picture}(140,160)(75,150)
\setlength{\unitlength}{2pt}   
\put(50,140){\makebox(20,6)[b]{F$_\sigma$-filter}}
\put(60,135){\vector(0,-1){10}} 
\put(50,120){\makebox(20,6)[b]{Strong P$^+$-filter}}
\put(60,115){\vector(-3,-1){30}}  
\put(60,115){\vector(3,-1){30}} 
\put(20,100){\makebox(20,6)[b]{Strong Baire Property}}
\put(80,100){\makebox(20,6)[b]{Completely Meager P$^+$-filter}} 
\put(30,95){\vector(3,-1){30}} 
\put(90,95){\vector(-3,-1){30}} 
\put(50,80){\makebox(20,6)[b]{Hereditarily Meager}}
\put(60,75){\vector(-3,-1){30}} 
\put(60,75){\vector(3,-1){30}}  
\put(20,60){\makebox(20,6)[b]{Completely Meager}}
\put(80,60){\makebox(20,6)[b]{Weakly Hereditarily Meager}}
\put(30,55){\vector(3,-1){30}} 
\put(90,55){\vector(-3,-1){30}}
\put(50,40){\makebox(20,6)[b]{Meager}}             
\put(140,60){\makebox(20,6)[b]{P$^+$-filter}}
\put(90,95){\vector(2,-1){59}}  
\put(150,55){\vector(-2,-1){59}}        
\put(60,35){\vector(3,-1){30}} 
\put(80,20){\makebox(20,6)[b]{Meager}}
\end{picture}

\vspace{4in}

The following lemmas will only reformulate these definitions into more
uniform and manageable ones; the proofs of the first three are left to the reader. 

\begin{lem} The following are equivalent for a given filter $\CF$:
\begin{enumerate}
\item $\CF$ is weakly hereditarily  meager.
\item For each sequence $\la X_i: i \in \omega \ra$ such that $\bigcup_i X_i \in \CF$
and each $X_i^c \in \CF$, 
\[ ( \exists n_0 < n_1< \cdots )(\forall Y \in \CF )( \fai k)\;
           Y \cap \bigcup_{i \in [n_k,n_{k+1})}X_i \neq \emptyset . \]
\end{enumerate}
\end{lem}

\begin{lem} The following are equivalent for a given filter $\CF$:
\begin{enumerate}
\item $\CF$ is completely  meager.
\item For each sequence of finite sets $\la X_i: i \in \omega \ra$
 such that $\bigcup_i X_i \in \CF^+$,
\[ ( \exists n_0 < n_1< \cdots )(\forall Y \in \CF )( \fai k)\;
           Y \cap \bigcup_{i \in [n_k,n_{k+1})}X_i \neq \emptyset . \]
\end{enumerate}
\end{lem}

\begin{lem} The following are equivalent for a given filter $\CF$:

\begin{enumerate}
\item $\CF$ is hereditarily meager.
\item For each sequence $\la X_i: i \in \omega \ra$ such that $\bigcup_i X_i \in \CF^+$
and each $X_i^c \in \CF$,
\[ ( \exists n_0 < n_1< \cdots )(\forall Y \in \CF )( \fai k)\;
           Y \cap \bigcup_{i \in [n_k,n_{k+1})}X_i \neq \emptyset. \]
\item For each sequence $\la X_i: i \in \omega \ra$ such that $\bigcup_i X_i
\in \CF$, $\bigcup_{i<k}X_i \notin \CF$ for each $k$ and only finitely
many $X_i$ belong to $\CF^+$, 
\[ (\exists n_0 < n_1< \cdots )(\forall Y \in \CF )( \fai k)\;
           Y \cap \bigcup_{i \in [n_k,n_{k+1})}X_i \neq \emptyset .\]

\end{enumerate}
\end{lem}

\begin{lem} The following are equivalent for a given filter $\CF$:

\begin{enumerate}
\item $\CF$ is a completely meager P$^+$-filter.
\item For each sequence $\la X_i: i \in \omega \ra$ such that $\bigcup_{i \geq n}
 X_i \in \CF^+$ for each $n$,
\[ (\exists h \in \fo)( \exists n_0 < n_1< \cdots )(\forall Y \in \CF )( \fai k)\;
           Y \cap \bigcup_{i \in [n_k,n_{k+1})}X_i \cap h(i) \neq \emptyset. \]
\item  For each sequence $\la X_i: i \in \omega \ra \sub \CF^+$, 
\[ (\exists h \in \fo)( \exists n_0 < n_1< \cdots )(\forall Y \in \CF )( \fai k)\;
           Y \cap \bigcup_{i \in [n_k,n_{k+1})}X_i \cap h(i) \neq \emptyset. \]
	   \
\end{enumerate}
\end{lem} 

\proof {\bf 1 $\implies$ 2:} Let $\CF$ be a completely meager P$^+$-filter
and fix a sequence $\la X_i: i \in \omega \ra$ such that $\bigcup_{i \geq n} X_i
\in \CF^+$ for each $n$. 

\noi Define $Y_n = \bigcup_{i \geq n} X_i$.  Then $Y_0 \sups Y_1 \sups
\cdots $ is a decreasing sequence from $\CF^+$ and therefore there is an
$X \in \CF^+$ almost included in each $Y_i$.  Since $\la \CF,X \ra$ must
be meager, choose $m_0 < m_1 < \cdots $ such that \[ ( \forall Y \in
\CF)( \fai k)\; Y \cap X \cap [m_k,m_{k+1}) \neq \emptyset .\]

\noi Now choose $n_0 < n_1 < \cdots$ and $h \in \fo$ such that
\[ (\forall k)(\exists \ell )\; [m_{\ell},m_{\ell +1}) \cap X \sub
         \bigcup_{i \in [n_k,n_{k+1})} X_i \cap h(i) .\]

\noi {\bf 2 $\implies$ 3} is immediate.

\noi {\bf 3 $\implies$ 1:} To show that $\CF$ is a P$^+$-filter, fix a
decreasing sequence $X_0 \sups X_1 \sups \cdots $ from $\CF^+$.
By assumption, there are $h \in \fo$ and a sequence $n_0 < n_1
< \cdots$ such that \[ (\forall Y \in \CF) \;
      Y \cap \bigcup_{i \in [n_k,n_{k+1})}X_i \cap h(i) \neq \emptyset. \] 

\noi If we let $X = \bigcup_k X_{n_k} \cap h(n_k) = \bigcup_i X_i \cap h(i)$, then 
$X \in \CF^+$ and is almost included in each $X_i$. 

\noi To show that $\CF$ is completely meager, let $X \in \CF^+$ be given and
use the assumptions on $X_i=X \sm i$. \qed

\vspace{.5in}

Now we are ready to prove the diagram implications; in view of the above
lemmas, all are straightforward except possibly two of them which we
prove in detail.

\begin{prop}
A strong P$^+$-filter has the strong Baire property.
\end{prop}

\proof We actually show that a strong P$^+$-filter has the Baire property
relative to any set $\CX \sub \ff$. Indeed, let 
\[ \CT = \{ s \in \fff: \CO_s \cap \CX \sub \CF \},\]

\noi and put $\CO = \bigcup_{s \in \CT} \CO_s$. We show that 
$(\CF \cap \CX) \triangle (\CO \cap \CX)$ is meager in $\CX$. 

\noi For this, let $\CS = \{s \in \fff: \CO_s \cap \CX \neq \emptyset$
and $\CO_s \cap \CX \nsubseteq \CF \}$; enumerate $\CS$ as $\la s_n :
n \in \omega \ra$. For each $n$, choose $A_n \in (\CO_{s_n} \cap \CX)
\sm \CF$ and put $X_n = A^c_n \sm lh(s) \in \CF^+$. By assumption,
there is an $h \in \fo$ such that for all $Y \in \CF$, $Y \cap X_n
\cap h(n) \neq \emptyset$ for all but finitely many $n$. Define now
\[M_k=\{X: (\fa \ell \geq k) \; X \cap X_\ell \neq \emptyset \} \sm
\CO . \]

\noi Clearly $(\CF \sm \CO) \sub \bigcup_k \CM_k$ and each $\CM_k$ is nowhere dense
in $\CX$. \qed

\begin{prop}
Every F$_\sigma$-filter is a strong P$^+$-filter.
\end{prop}

\proof Let $\CF=\bigcup_n \CC_n$ be a countable increasing union of closed sets
and $\la X_n : n \in \omega \ra \sub \CF^+$. Put $\CC= \bigcup \{X
\cup n : X \in \CC_n \}$ and define $\CC'_n= \CC \cup X_n$. Choose $h
\in \fo$ such that for all $n$,
 \[ (\fa Y,Z \in \CC'_n) \; Y \cap Z \cap \lsb n,h(n) \rb \neq \emptyset . \]

\noi The existence of such a function is guaranteed by the fact that $\CC'_n$
is closed. Now if $Y \in \CF$, say $Y \in \CC_i$, then $Y \cup i \in \CC$ 
and therefore $(Y \cup i) \cap X_n \cap h(n) \neq \emptyset$ for all $n$,
and thus $Y \cap X_n \cap h(n) \neq \emptyset$ for all $n \geq i$ as desired. \qed

\section{Examples}

We now produce examples showing that many implications in the diagram
are not reversible.

\begin{eple}
A weakly hereditarily meager filter which is not completely meager.
\end{eple}

Let $\CU$ be any ultrafilter not containing say the even numbers $E$. Define 
$\CF= \{E \cup X: X \in \CU \}$. Then $\la \CF, E^c \ra = \CU$ is not meager
and thus $\CF$ is not completely meager.

\noi Now let $\la X_i:i \in \omega \ra$ be given such that 
$\bigcup_i X_i \in \CF$ and each $X_i^c \in \CF$. Choose a sequence
$n_0 < n_1 \cdots $ such that for each $k$,
$(E \sm k) \cap \bigcup_{i \in [n_k,n_{k+1})} X_i \neq \emptyset$.
This is possible as $E \cap \bigcup_{i<k}X_i$ is finite for each $k$.
Since any $Y \in \CF$ almost contains $E$, we conclude that $\CF$ is 
weakly hereditarily meager by lemma 2.2. \qed

\begin{eple}
A completely meager filter which is not weakly hereditarily meager.
\end{eple}

This example comes essentially from \cite{tal2}, where a completely meager filter
without the strong Baire property is produced.
Let again $\CU$ be an ultrafilter and let $\CF = \CU \otimes {\sl Fr}$. 

\noi {\bf Claim 1:} $\CF$ is completely meager.

Indeed let $X \sub \omega \times \omega$ be compatible with $\CF$; this
means that $X_n=\{(n,m):m \in X\}$ is infinite for infinitely many $n$
(actually on a set in $\CU$).  Define a sequence of finite sets $\la
s_i:i \in \omega \ra$ such that $(s_{i+1} \sm s_i) \cap X_j \neq \emptyset$
for each $j\leq i$ such that $X_j$ is infinite. Clearly

\[ (\forall Y \in \CF)(\fai i) \; Y \cap X \cap s_i \neq \emptyset \]
\noi and thus $\la \CF, X \ra$ is meager.

\noi {\bf Claim 2:} $\CF$ is not weakly hereditarily meager.

Define $f \in \fo$ by $f(n,m)=n$, the downward projection. Then no
$f^{-1}(i)$ is compatible with $\CF$ but  $f(\CF)=\CU$ is not meager.
\qed

\begin{eple}
A completely meager P$^+$-filter which is not a strong P$^+$-filter.
\end{eple}

We build a mad family $\CA$ such that the filter $\CF(\CA)$ has the
required properties.  It is essentially proved in \cite{laf} that any
filter of the form $\CF(\CA)$ will necessarily be a completely meager
P$^+$-filter; we include the proof for completeness. 

\noi {\bf CLAIM 1:} All filters of the form $\CF(\CA)$ are competely
meager P$^+$-filters.

Notice that a set $X \in \CF(\CA)^+$ if and only if $X$ has infinite
intersection with infinitely many members of $\CA$. Now if $Y_0 \sups
Y_1 \sups \cdots$ are from $\CF(\CA)^+$, choose $\{A_n: n \in \omega
\} \sub \CA$ such $Y_n \cap A_m \neq \emptyset$ for each $n,m$;
finally define a sequence $\la n_k:k \in \omega \ra$ such that \[ Y_k
\cap A_\ell \cap [n_k,n_{k+1}) \neq \emptyset \] \noi for each $k$ and
$\ell \leq k$. Clearly $Y= \bigcup_k Y_k \cap [n_k,n_{k+1}) \in
\CF(\CA)^+$ and is almost included in each $Y_n$. Actually we see
that for any $X \in \CF(\CA)$, $X \cap Y \cap [n_k,n_{k+1}) \neq \emptyset$
for all but finitely many $k$ and thus we have also proved that
$\CF(\CA)$ is completely meager.

Now to make sure that $\CF(\CA)$ is not a strong P$^+$-filter, first
fix an almost disjoint family $\CB=\{B_\alpha: \alpha < \gc \}$ of
size continuum and fix an enumeration of $\fo=\{f_\alpha:  \alpha <
$\gc$ \}$. Fix further a partition  $\la X_n : n \in \omega
\ra$ of $\omega$. Now we define an almost disjoint family 
\[ \CA'=\{
\bigcup_{n \in B_\alpha} X_n \cap f_\alpha (n) : \alpha < \mbox{\gc} \}. \]

\noi Let $\CA$ be any mad family extending $\CA'$ together with a
partition of the $X_n$'s into countably many infinite sets (to make
sure that $X_n \in \CF(\CA)^+$). Clearly $\CF(\CA)$ is not a strong
P$^+$-filter.

\begin{eple}
A filter with the strong Baire property which is not a P$^+$-filter.
\end{eple}

Let $\CF={\sl Fr} \otimes {\sl Fr}$. Then $\CF$ is a Borel filter
and thus has the strong Baire property. If we define $Y_n =
\{(\ell,m): \ell \geq n \} \in \CF$, then $\CF$ contains the
complement of all $Y$ almost included in each $Y_n$.

\begin{eple}
A non-meager P$^+$-filter.
\end{eple}

It is worth noticing that the concept of P$^+$-filter and the traditional
notion of P-filter are incomparable, that is one does not imply the
other.  The difficult question of producing a non-meager P-filter has
been around for a some time, see \cite{just} for some results  and
applications in this area.

We show that Simon's example (\cite{simon}) of a non-meager filter
generated by \gb sets is actually a P$^+$-filter.  As in his
construction, we consider two possibilities. 

\noi {\bf Case 1: (\gb $<$ \gd)} Choose an unbounded family $\la
f_\alpha: \alpha < $ \gb $ \ra$ of strictly increasing functions such that
\[ \alpha < \beta \implies f_\alpha <^* f_\beta. \]

\noi Fix now $g \in \fo$ not dominated by any $f_\alpha$, put
\[X_\alpha = \{n: f_\alpha(n) < g(n) \} \] \noi and define $\CF =\la
X_\alpha: \alpha < $ \gb $ \ra$. $\CF$ is a P-filter and the next lemma
shows that it really is a P$^+$-filter.

\begin{lem}
Any filter generated by less that \gd  sets is a P$^+$-filter.
\end{lem}

Let $\CH$ be a filter generated by $\la X_\alpha :\alpha < \lambda \ra$
for some $\lambda < $ \gd,  and fix a descending sequence $Y_0 \sups Y_1
\sups \cdots $ from $\CH^+$. For each $\alpha < \lambda$, choose $h_\alpha
\in \fo$ such that
\[ (\fa n) (Y_n \cap X_\alpha \cap \lsb n,h_\alpha(n) \rb \neq \emptyset. \]

\noi Since $\lambda <$\gd, fix $h \in \fo$ not dominated by any
$h_\alpha$ and let $Y=\bigcup_n Y_n \cap h(n)$ almost included in each
$Y_n$. But for each $\alpha$, if $h_\alpha(n) \leq h(n)$, $Y \cap
X_\alpha \cap \lsb n,h(n) \rb \neq \emptyset$ and therefore $Y \in \CH^+$. \qed

Finally we must show that $\CF$ is non-meager. So suppose we are given
a sequence $n_0 < n_1 < \cdots $ such that 

\[ (\fa X \in \CF)(\fai k) \; X \cap [n_k,n_{k+1}) \neq \emptyset .\] 

\noi Then for each $\alpha$, choose $m$ such that $X_\alpha \cap
[n_k,n_{k+1}) \neq \emptyset$ for each $k \geq m$. Therefore for each
$\ell$, $f_\alpha(\ell) \leq n_{\ell + m +1}$ and thus for all $\ell
\geq m+1$ we have $f_\alpha(\ell) \leq n_{2 \ell}$. We conclude that
the family  $\la f_\alpha: \alpha < $\gb $ \ra$ is bounded by
$f(\ell)=n_{2 \ell}$, a contradiction.

\noi {\bf Case 2: (\gb =\gd)} Choose a scale  $\la f_\alpha: \alpha <
$\gb=\gd$ \ra$ such that 

\begin{enumerate}

\item Each $f_\alpha$ is strictly increasing and positive.
\item $(\fa \beta < \alpha)(\fai n)(\ex \ell)
              f_\alpha(n) < f_\beta(\ell) < f_\beta(\ell +1) < f_\alpha(n+1)$.

\end{enumerate}

\noi For each $\alpha$, put $X_\alpha= \bigcup_n \lsb
f^{2n}_\alpha(0),f^{2n+1}_\alpha(0) \rb $ and define $\CF=\la X_\alpha :
\alpha < $\gb=\gd $\ra$. We must show that $\CF$ is a non-meager
P$^+$-filter.

\begin{lem}
$(\fa A \in \so)(\ex \alpha)(\fai n) \; A \cap
\lsb f^n_\alpha(0),f^{n+1}_\alpha(0) \rb \neq \emptyset.$
\end{lem}

\proof Given $A$, define $f \in \fo$ such that $A \cap \lsb n,f(n) \rb
\neq \emptyset$ for each $n$. Now there must be an $\alpha$ such that
$f \leq^* f_\alpha$ and therefore
\[ (\fai n) \; A \cap \lsb n,f_\alpha(n) \rb \neq \emptyset .\]

\noi Thus  $(\fai n) \; A \cap \lsb f_\alpha
(n),f_\alpha(f_\alpha(n)) \rb \neq \emptyset $, and finally 
\[ (\fai n)\; A \cap \lsb f^n_\alpha(0),f^{n+1}_\alpha(0) \rb \neq \emptyset \]

\noi as desired. \qed

We are ready to show that $\CF$ is a P$^+$-filter. Fix a descending
sequence $Y_0 \sups Y_1 \sups \cdots$ from $\CF^+$. By induction,
define a sequence of ordinals $\alpha_n < $\gb=\gd by $\alpha_0=0$ and
given $\alpha_n$, let $\CF_n$ be the filter generated by \[ \{X_\gamma
: \gamma < \alpha_n \} \cup \{ Y_k: k \in \omega \}. \]

\noi Now choose $\alpha_{n+1} > \alpha_n$ so that
\[ (\fa Y \in \CF_n)(\fai k) \; Y \cap 
 \lsb f^k_{\alpha_{n+1}}(0), f^{k+1}_{\alpha_{n+1}}(0)\rb \neq
    \emptyset .\]

\noi This is possible using the lemma and the fact that $\CF_n$ is
generated by less than \gb=\gd  sets. Finally let $\alpha= \sup
\{\alpha_n: n \in \omega \}$.  We thus get

\[ (\fa Y \in \CF_\alpha :=
\la X_\gamma: \gamma < \alpha \ra \cup \la Y_n:n \in \omega \ra)
    (\fai k) \; Y \cap \lsb f^k_{\alpha}(0), f^{k+1}_{\alpha}(0)\rb
           \neq    \emptyset .\]

\noi Therefore it suffices to find $Z \sub^* Y_n$ such that

\begin{enumerate}
\item $Z \in \CF^+_\alpha$
\item $(\fai k)\; 
           Z \cap  \lsb f^k_{\alpha}(0), f^{k+1}_{\alpha}(0)\rb
           \neq    \emptyset.$
\end{enumerate}

\noi But again, $\CF_\alpha$ is generated by less than \gb sets, say
by $\la X_\beta: \beta < \lambda \ra$ for some $\lambda <$ \gb, and we may
assume that this collection is also closed under finite intersections.
Now for each such $\beta$, choose $h_\beta \in \fo$ such that $Z_\beta
= \bigcup_n Y_n \cap h_\beta(n)$ satisfies: 

\begin{enumerate}

\item $(\fai k) \; Z_\beta \cap \lsb f^k_\beta(0), f^{k+1}_\beta(0) \rb \neq
       \emptyset.$ 
\item $Z_\beta \cap X_\beta$ is infinite.   

\end{enumerate}

\noi Now choose $h \;^*\!>h_\beta$ for $\beta < \lambda$. Then
$Z=\bigcup_n Y_n \cap h(n)$ is as desired.

\begin{eple}
(CH) A strong P$^+$-filter not included in any F$_\sigma$ filter.
\end{eple}

Enumerate all closed sets which generate a proper filter as $\la
\CC_\alpha: \alpha < \omega_1 \ra$ and all sequnces of potential
candidates for members of $\CF^+$ as $\la X^\alpha_n:n \in \omega \ra$
for $\alpha < \omega_1$. We build a filter $\CF$ in $\omega_1$ stages.
At stage $\alpha$, assume that we have built for $\beta < \alpha$: 

\begin{enumerate} 
\item $\la a^\beta_n:n \in \omega \ra$ a sequence of finite sets. 
\item $X_\beta \in \so$ such that $\lim_n |X_\beta \cap a^\gamma_n| = + \infty$
for all $\gamma < \alpha$.
\item $\beta < \gamma < \alpha \implies X_\gamma \sub^* X_\beta$.
\end{enumerate}

\noi Let $\CF_\alpha$ be the filter generated by $\{X^c: (\fa n) |X
\cap a^\beta_n| \leq m \}$ for $m \in \omega$ and by $\{X_\beta: \beta
< \alpha \}$. We are given a closed set $\CC = \CC_\alpha$ and  we
want to make sure that $\CF_\alpha \nsubseteq \la \CC \ra$. So we
suppose that in fact that $\CF_\alpha \subseteq \la \CC \ra$ and we
may as well assume that $\alpha = \omega$.

\noi Choose a sequence $n_0 < n_1 < \cdots $ such that 
\[ ( \fa Y_1,\dots ,Y_k \in \CC ) \; Y_1 \cap \cdots Y_k \cap [n_k,n_{k+1}) \neq \emptyset.\]

\noi {\bf Claim:} $(\fa \beta < \alpha)(\fa m)(\fa Z)$, if $|Z \cap a^\beta_k| \leq m$
for each $k$, then 
\[ (\fa Y \in \CC)(\fai \ell)\; (Y \cap [n_\ell,n_{\ell + 1})) \sm Z \neq \emptyset. \]

\proof Indeed we have made sure to put $Z^c$ in $\CF_\alpha$ for all
such $Z$'s, and $\CF_\alpha \sub \la \CC \ra$ by assumptions.  \qed

\noi The construction of $X_\alpha$ can therefore be obtained by a
dovetailed construction.  Now to handle the sequence  $\la
X^\alpha_n:n \in \omega \ra $, given that each $X^\alpha_n \in \la
\CF_\alpha, X_\alpha \ra^+$, choose $f \in \fo$ such that

\[ (\fa n) \; |X_\alpha \cap X^\alpha_n \cap \lsb n,f(n) \rb | \geq n, \] 

\noi and define $a^\alpha_n = \lsb n,f(n) \rb \cap X^\alpha_n$.

This completes the construction. $\CF = \bigcup_\alpha \CF_\alpha$ is the desired filter.

\begin{eple}
(CH) A mad family $\CA$ such that $\CF(\CA)$ is a strong P$^+$-filter.
\end{eple}

In the following, $\CI(\CA)$ denotes the ideal generated by members of $\CA$.
We start with a lemma.

\begin{lem}
(CH) There is a mad family $\CA$ such that
\[ (\fa X \notin \CF(\CA) \cup \CI(\CA))(\ex \alpha)(\fa \beta > \alpha)\;
          A_\beta \sm X \mbox{ and } A_\beta \cap X \mbox{ are infinite }. \]
\end{lem}

\proof Start with a partition of $\omega$ into infinite sets $\la
A_n:n \in \omega \ra$.  Now enumerate $\so$ as $\la X_\alpha: \alpha <
\omega_1 \ra$. Suppose that at stage $\alpha <\omega_1$, we have
$\CA_\alpha = \la A_\beta: \beta < \alpha \ra$ and $\CX_\alpha \sub
\so \sm [\CF(\CA_\alpha) \cup \CI(\CA_\alpha)]$. Relist $\CA_\alpha$
as $\la B_n:n \in \omega \ra$.

\noi {\bf Case 1:} $X_\alpha \in \CF(\CA_\alpha) \cup \CI(\CA_\alpha)$.
There is nothing to do.

\noi {\bf Case 2:} $X_\alpha \notin \CF(\CA_\alpha) \cup
\CI(\CA_\alpha)$.  Thus in particular both $X_\alpha \nsubseteq^* B_0
\cup \cdots B_n$ and    $X^c_\alpha \nsubseteq^* B_0 \cup \cdots B_n$
for each $n$.

\noi List $\CX_\alpha \cup \{X_\alpha\}$ as $\la X_n:n \in \omega \ra$
and we build $A_\alpha$ such that:

\begin{enumerate}
\item $(\fa n) \; A_\alpha \cap B_n =^* \emptyset$.
\item $(\fa n) \; A_\alpha \sm X_n \mbox{ and } A_\alpha \cap X_n \mbox{ are infinite }.$
\end{enumerate}

\noi At stage $n$, assume that we have $a_n \in \fso$, then choose
\[ k_i \in X_i \sm (B_0 \cup \cdots B_n \cup n) \mbox{ for }i<n,\]
\[ \ell_i \in X^c_i \sm (B_0 \cup \cdots B_n \cup n) \mbox{ for }i<n,\]

\noi and let $a_{n+1}=a_n \cup \{k_i,\ell_i\}_{i<n}$. Hence $A_\alpha= \bigcup_n a_n$ satisfies
1-2 and put $\CA_{\alpha+1}=\CA_\alpha \cup \{A_\alpha\}$, $\CX_{\alpha + 1}=
\CX_\alpha \cup \{X_\alpha \} \sub \so \sm [\CF(\CA_{\alpha+1}) \cup \CI(\CA_{\alpha+1})]$.
This completes the construction. The verification that $\CA=\{A_\alpha: \alpha < \omega_1 \}$
has the desired properties is straightforward. \qed

Now we show that such a mad family yields a strong P$^+$-filter
$\CF(\CA)$.  But given a sequence $\la X_n: n \in \omega \ra$ from
$\CF(\CA)^+$, we can choose by the assumptions on $\CA$ countably
members $\la A_n : n \in \omega \ra$ from $\CA$ such that $A_n \cap X_m$
is infinite for each $n$ and $m$.  Now define $f \in \fo$ such that 

\[(\fa n)(\fa i<n) \; A_i \cap X_n \cap \lsb n,f(n) \rb \neq \emptyset. \]

\noi If $Y \in \CF(\CA)$, there must be an $A_n$ such that $A_n \sub^* Y$ and therefore
$Y \cap X_n \cap f(n) \neq \emptyset$ for all but finitely many $n$.

\section{Conclusion}

The two implications that remain to be solved are the following.

\noi {\bf Question 1:} Does a completely meager P$^+$-filter, or
actually even a hereditarily meager filter necesarily have the strong
Baire property?

One remains in ZFC.

\noi {\bf Question 2:} Is there a strong P$^+$-filter not included in
any F$_\sigma$-filter?

We still have very little understanding of the structure of filters of the
form $\CF(\CA)$ for a mad family $\CA$. In particular, the following remains 
unsolved.

\noi {\bf Question 3:} Do all filters of the form $\CF(\CA)$ for a mad family
$\CA$ have the strong Baire property?

We however have a meagerness property that separates filters of the form
$\CF(\CA)$ and F$_\sigma$-filters, namely the strong P$^+$ property.  In
view of the results of \cite{laf} that F$_\sigma$-filters can be
diagonalized in a forcing extension without adding an unbounded real, we
may ask the following. 

\noi {\bf Question 4:} Can all strong P$^+$-filters be diagonalized in
a forcing extension without adding unbounded reals?

Of course the fundamental motivation for these questions is:

\noi {\bf Question 5:} Can filters of the form $\CF(\CA)$ for a mad family
$\CA$ be diagonalized in a forcing extension without adding an
unbounded real?

\end{document}